\documentclass[12pt]{article}         
\usepackage{amsfonts,amssymb} 
\makeatletter
\def\draft{\textheight=10.5truein \textwidth=7.5truein \parindent=8pt
           \voffset=-1truein \topmargin=0Truein
           \ifcase \@ptsize \hoffset=-1.5truein \or \hoffset=-1.35truein
                        \or \hoffset=-1.15truein \fi}
\def\quality{\textheight=240mm \textwidth=160mm \topmargin=0Truein
             \ifcase \@ptsize \hoffset=-23mm \or \hoffset=-20mm
                          \or \hoffset=-15mm \fi}
\makeatother
\quality %

\def\beq#1#2{\begin{equation} \label{#1} #2 \end{equation}}
\def\bea#1{\begin{eqnarray*} #1 \end{eqnarray*}} \def\a{\!\!\!&\!\!\!\!&}
\def\function#1{\left\{\!\!\!\begin{array}{ll} #1 \end{array} \right.}

\def\proof{\smallskip \noindent {\bf Proof. \ }}       
\def\blanksquare{\,\,\,$\sqcup\!\!\!\!\sqcap$}         
\def\qed{\hfill\blanksquare\linebreak\smallskip\par}   

\def\thname{Theorem}     \def\lmname{Lemma}      \def\prname{Proposition}
\def\dfname{Definition}  \def\crname{Corollary}  \def\rmname{Remark}

\newtheorem{theorem}{\thname}[section]   
\newtheorem{lemma}{\lmname}[section]     
\newtheorem{corollary}[lemma]{\crname}   

\newtheorem{dftn}{\dfname}[section]
\newtheorem{rmrk}[lemma]{\rmname}

             \catcode`@=11 \@addtoreset{equation}{section} \catcode`@=12

\def\bline(#1,#2)(#3,#4)(#5){\put(#1,#2){\line(#3,#4){#5}}}  
\newcommand\mlbscale{1pt} 
\newif\iffigs\figstrue 

\def\bfig(#1,#2)#3#4{\begin{figure} \begin{center}
    \framebox{\setlength{\unitlength}{\mlbscale}
       \iffigs \begin{picture}(#1,#2) #3 \end{picture}
       \else \begin{picture}(60,10)(0,0)
                   \put(0,0){\framebox(60,10){Figure}} \end{picture} \fi}
    \end{center} \caption{#4} \end{figure}}

\def\Bfig(#1,#2)#3#4{\begin{figure} \begin{center}
    \setlength{\unitlength}{\mlbscale}
       \iffigs \begin{picture}(#1,#2) #3 \end{picture}
       \else \begin{picture}(60,10)(0,0)
                   \put(0,0){\framebox(60,10){Figure}} \end{picture} \fi
    \end{center} \caption{#4} \end{figure}}

\def\bpic(#1,#2)#3{\setlength{\unitlength}{\mlbscale}
    \begin{picture}(#1,#2) #3 \end{picture}}

\def\IR{{\mathbb{R}}} \def\ep{\varepsilon}  \def\phi{\varphi}
\def\la{\lambda}      \def\map{T}           \def\IZ{{\mathbb{Z}}}
\def\cM{{\cal M}}     \def\n{\noindent}     \def\t{\tilde}
\def\IL{{\bf L}}      \def\IC{{\bf C}}
\def\toas#1{\stackrel{#1}{\longrightarrow}}
\def\v#1{\overline{#1}} 
\def\vmap{\v{T}} \def\vX{\v{X}} \def\vx{\v{x}} \def\vy{\v{y}}
\def\vm{\v{m}} \def\vmu{\v{\mu}} \def\vQ{\v{Q}} 
 \def\ga{\gamma} 

\def\cB{{\cal B}}  \def\Q{\v{Q}_\ep}
\def\cF{{\cal F}} \def\V{{\bf V}} \def\La{\Lambda}  \def\al{\alpha}
\def\*#1{#1^*}    \def\0#1{\breve#1}  \def\2#1{\acute#1}

\def\ns#1{||#1||} \def\nw#1{|#1|}

\def\?#1{} 

 {\vskip-0.5truecm}

\begin{document}
\title{Condensation versus independence \\ in weakly interacting CMLs}
\author{Michael Blank\thanks{
        Russian Academy of Sci., Inst. for
        Information Transm. Problems,
        and Laboratoire Cassiopee UMR6202, CNRS, ~
        e-mail: blank@iitp.ru}
        \thanks{This research has been partially supported
                by Russian Foundation for Fundamental Research
                and French Ministry of Education grants.}
       }
\date{25 February 2010} 
\maketitle

\begin{abstract}
We propose a simple model unifying two major approaches to the
analysis of large multicomponent systems: interacting particle
systems (IPS) and couple map lattices (CML) and show that in the
weak interaction limit depending on fine properties of the
interaction potential this model may demonstrate both
condensation/synchronization and independent motions. Note that
one of the main paradigms of the CML theory is that the latter
behavior is supposed to be generic. %
The model under consideration is related to dynamical networks and
sheds a new light to the problem of synchronization under weak
interactions. \end{abstract}

\section{Introduction}\label{s:intro}
At present there are two major approaches to the analysis of large
multicomponent systems: Interacting Particle Systems (IPS) and
Coupled Map Lattices (CML). Ideas and especially methods used in
these approaches are strikingly different and it is hard to find a
single article which discuss their connections on a reasonably
serious level. Nevertheless in this paper we present a simple
mathematical model enjoying interpretations in terms of both of
these approaches which helps a lot in its study.

Let $X$ be a compact subset of $\IR^d$ equipped with the Euclidian
metric $\rho(\cdot,\cdot)$ and let
$(\vmap,\vX):=\bigotimes_i^N(\map,X)$ be a direct product of $N$
copies (in general $0<N\le\infty$) of the measurable nonsingular
dynamical system $(\map,X,\cB,\mu_\map)$. Consider also a family
of measurable maps $\vQ_\ep:\vX\to\vX$ to which we shall refer as
{\em interactions}. The parameter $\ep\ge0$ here measures how far
the map $\vQ_\ep$ is from the
identical map, namely %
$\ep:=\sup_{\vx\in\vX}\v\rho(\vQ_\ep\vx,\vx)$, where
$\v\rho(\vx,\v{y}):=\max_i \rho(x_i,y_i)$.

Now a CML is defined as a composition of the direct product map
and the interaction: $\vmap_\ep:=\vmap\circ\vQ_\ep$ and the system
without interactions can be formally written as
$\vmap_0\equiv\vmap$.

A typical example of interactions considered in the literature is
the so called ``diffusive coupling'':
\beq{e:diffusive}{(\vQ_\ep\vx)_i:=\ep x_i +
  \frac{1-\ep}{|J_i|}~{\sum_{j\in J_i}x_j} ,} %
where the summation is taken over ``neighboring elements'' $J_i$
of the local unit $i$. The formula~(\ref{e:diffusive}) indicates
that speaking about a CML one usually has in mind a certain
structure or topology of connections, say a fixed ordered graph of
interactions (see e.g.
\cite{Bl-mon,Bl21,BB,BK,Bun-C,Bun-S,Bu,KL05,Wu}). Another
possibility is to consider a dynamic setup allowing a dynamical
switching of connections. Some results in this direction were
ob6ained in \cite{LAJ} (in the case of time-varying couplings) and
in \cite{KL08}. In what follows we propose a simple model of a
dynamical switching and show that in the weak interaction limit
depending on fine properties of the interaction potential this
model may demonstrate both condensation/synchronization and
independent motions.

Motivated by the idea of the interaction between the local systems
induced by (somewhat artificial)\footnote{The interaction
    discussed in \cite{KL08} resembles more closely the exchange of
    velocities between colliding particles rather than the collision
    itself. The latter is modelled by a prior arrangement of
    non-overlapping ``traps'' for pairs of particles.} %
``collisions'' introduced in \cite{KL08} we interpret the CML as
follows. A {\em configuration} $\vmap_\ep^t\vx:=\vx^t=\{x_i^t\}$
at time $t\ge0$ is associated to the collection of identical
particles with unit masses located at points $x_i^t\in X$ which
wander in the common domain $X$ independently until they come
close enough to each other. As an example one can think about a
billiard type system. The interaction occurs only when the
particles are coming $\ep$-close to each other and consists in the
attraction to the local common center of gravity.

Assume that $(X,\cB,m,\rho)$ is a compact convex Lebesgue metric
space $X\in\IR^d$. For a configuration $\vx\in\vX$ and each index
$i$ denote by $J_i=J_i(\vx):=\{j:~~\rho(x_i,x_j)\le\ep\}$ the set
of indices of particles ``interacting'' with the $i$-th
one.\footnote{Later we shall consider some other choices of the
     set $J_i$.} %
Then we define the $i$-th coordinate of the interaction map as
follows: %
\beq{e:interaction}{(\vQ_\ep\vx)_i:=\ga x_i +
  \frac{1-\ga}{|J_i|}~{\sum_{j\in J_i}x_j} .} %
The parameter $\ga\in[0,1]$ and if $\ga=0$ only the second
summand, corresponding the ``center of gravity'' of particles
belonging to the $\ep$-neighborhood of the point $x_i$, survives
in the expression above. Therefore the physical meaning of the
interaction under study is that the $i$-th particle is moved in
the direction of the center of gravity of particles belonging to
the $\ep$-neighborhood of the $i$-th particle.


If $\ga\in(0,1)$ we say that the interaction is {\em soft}, $\ga=0$
corresponds to the situation when a particle jumps  directly to
the center of gravity, which we call {\em rigid} by contrast.
Finally $\ga=1$ nullifies the interaction.

\?{In general one may consider an arbitrary interaction potential
$U(\cdot)$ instead of the indicator function, but we restrict the
analysis to this special case in particular due to its evident
analogy to the celestial mechanics setting.}

Our main results may be formulated as follows.
Let the $(\map,X,\cB,\mu_\map)$ be a weakly mixing measurable dynamical
system and let $\mu_\map$ be its only Sinai-Bowen-Ruelle (SBR) measure
(see definitions and discussion in Section~\ref{s:defs}).
These assumptions are enough for the case of rigid interactions,
but to study soft interactions one inevitably needs some additional smoothness type assumptions for the local map (see Section~\ref{s:soft1}).

\begin{theorem}\label{t:main1} \label{t:sbr2}
For each $N\in\IZ_+$ there are constants $\ep_\map>0, ~\ga_\map>0$
such that $\forall~ 0\le\ep<\ep_\map,~ 0\le\ga<\ga_\map$ the CML
$(\vmap_\ep,\vX)$ has a SBR measure $\vmu$ (which does not depend
on $\ep$), is supported by the ``diagonal'' of $\vX$, and whose
projections to each ``coordinate'' coincide with $\mu_\map$. If $N=2$
this SBR measure is unique.
\end{theorem}%

This result demonstrates that arbitrary weak interactions of
a CML with a ``generic'' local map may demonstrate discontinuity
at $0$ of the SBR measure $\vmu_\ep$ considered as a function on $\ep$.
Previously such results were known only for ``wild'' examples of
local maps having periodic turning points type singularities, see e.g. \cite{Bl21,Bl-mon}.
From the point of view of interacting particle systems this
theorem can be interpreted as a kind of a condensation phenomenon
when particles are gathering together under dynamics.

In the opposite case when $\ga$ is close to $1$ we observe
a more ``classical'' effect, namely under some reasonably general
technical assumptions specified in Section~\ref{s:soft2}
the only SBR measure of the CML converges weakly to the direct
product measure as $\ep\to0$. Bassicaly we need that the action
of the transfer operator corresponding to the local map in a
suitable Banach space of signed measures be quasi-compact, i.e.
can be decomposed into a sum of a $\theta$-contracting operator
(with $\theta<1$) and a finite dimensional projector.

\begin{theorem}\label{t:main2} For each $N\in\IZ_+$ there is
$\ep'_\map>0$ such that $\forall~ 0\le\ep<\ep'_\map$ and each
$\ga\in(2^N\theta,1]$ the CML $(\vmap_\ep,\vX)$ has a unique SBR
measure
$\vmu_{\ep,\ga}\toas{\ep\to0}\vmu_\map
:=\underbrace{\mu_\map\otimes\dots\otimes\mu_\map}_N$.
\end{theorem}%

From the interacting particles interpretation this result tells
that in this case weak interactions lead to independent particles
motions. Thus soft interactions demonstrate very different and
more ``classical'' effects in comparison to the rigid
interactions. Note that here one may consider not only particles
of different masses but even the local maps need not be identical.

To understand better the nature of the phenomena under study
let us consider the behavior of the CML in more detail and under
a bit different point of view which is more closely related to the
IPS approach.

We say that a {\em synchronization} or {\em condensation} with the
{\em basin} $\v{Y}$ takes place in the CML $(\vmap_\ep,\vX)$ if %
$\limsup\limits_{t\to\infty}\max\limits_{i,j}\rho(x_i^t,x_j^t)=0$ %
for each configuration $\vx\in\v{Y}$.
Similarly a {\em desynchronization} with the {\em basin}
$\v{Y}$ means that %
$\liminf\limits_{t\to\infty}\max\limits_{i,j}\rho(x_i^t,x_j^t)>0$ %
for each $\vx\in\v{Y}$.
We say also that a certain property with the domain $A$ and a
probabilistic measure $\nu$ is $\nu$-{\em global} if $\nu(A)=1$.

In words, the synchronization means that the dynamics converges to
the {\em diagonal} $\v{D}:=\{\vx\in\vX:~~x_i=x_j~\forall i,j\}$,
while the desynchronization corresponds to the absence of such
convergence. We refer the reader to \cite{PRK} for the discussion
of various physical aspects of the synchronization.

Let us start with the rigid case. The following statement
formulated in terms of the de/synchronization not only clarifies
the situation but also gives some additional insight about the
dynamics.

\begin{theorem}\label{t:sync2} Under the assumptions of
Theorem~\ref{t:sbr2} $\forall\ep>0$ the synchronization with the
open basin $\v{Y}$ of positive product measure $\vmu_\map$ takes
place. Moreover for $N=2$ this synchronization is
$\v\mu_\map$-global. On the other hand, for $N\ge3$ there is an
analytic local map $\map$ for which the desynchronization with an
open positive product measure $\vmu_\map$ basin takes place.
\end{theorem}

It is worth note that our numerical simulations show that for all
1D and 2D mixing maps that we tried for {\em all} initial particle
configurations the global synchronization was observed. This is
especially striking because first one might expect the
synchronization only for almost all initial particle
configurations, and second even for the cases where we are able to
prove the desynchronization we do not see it in the numerical
experiment. The explanation is that arbitrary small round-off
errors may change the behavior of a chaotic dynamical system
drastically (see e.g. \cite{Bl-mon}) and they are responsible in
this case as well.

The proof of Theorem~\ref{t:sync2} is based on the following
observations. If two particles are coming $\ep$-close to each
other, then after the interaction their positions coincide.
Therefore if $N=2$ we only need to show that for almost all
2-particle configurations their trajectories will hit
simultaneously the $\ep$-neighborhood of the diagonal $\v{D}$. To
prove the latter statement one uses that outside of this
$\ep$-neighborhood there are no interactions and hence it is
enough to check the same statement for the direct product system
$(\vmap,\vX,\v\cB,\vmu_\map)$, which in turn follows from the weak
mixing of the local map.

One is tempted to extend this construction directly for the
case $N>2$ along the following lines. %
First, if the trajectory of an initial configuration hits an
$\ep$-neighborhood of the diagonal $\v{D}$ (as in the case $N=2$)
then after the interaction the synchronization takes place: all
particles will share the same position. %
Second, additionally to the previous argument one considers
multiple (triple, etc.) collisions between particles (which occur
near secondary diagonals of $\vX$) expecting that each collision
reduces the number of unmatched particles.

Unfortunately, there are two fundamental obstacles to this naive
approach. First, a single collision of $n>2$ particles does not
necessarily imply that they will share the same location after the
collision. Indeed, assume that we have $n>2$ particles uniformly
distributed along a circle of radius $r_n=\frac\ep2/\sin(\pi/n)$.
Then in the $\ep$-neighborhood of each particle there are two
neighbors located at distance $\ep$. Therefore after the
interaction instead of coming to a single common center of
gravity, the particles will be again uniformly distributed along a
circle of only a slightly smaller (provided $n\gg1$) radius
$r_n':=r(1+2\cos(2\pi/n))/3$, e.g. $r'_4=r_4/3$ and $r'_5\approx
r_5\times0.539$. Observe that we do not assume any smoothness of
the local map $\map$, and thus after the next application of
$\map$ neighboring points might become arbitrary far from each
other.

The second and even more important obstacle is that even a single
collision changes the trajectory of a particle and one cannot
apply (at least directly) arguments related to ergodic properties
of the original local system $(\map,X,\cB,\mu)$. Moreover after each
decrease of the number of unmatched particles one needs to check
that the system is still in a ``general position'' in order to use
again the mixing property.

The result about the desynchronization shows that the smoothness
of the map $\map$ does not cure these pathologies even in the
simplest setting. In a sense the doubling map turns out to be the
``worst'' local one-dimensional map for our problem.

\bigskip

Now we turn to the analysis of soft interactions. The most striking
difference to the rigid case is that results of Theorem~\ref{t:sbr2}
or \ref{t:sync2} are no longer available without some assumptions
on the smoothness of the local system. Nevertheless
Theorem~\ref{t:soft-rigid} (Section~\ref{s:soft1}) shows that if
the local map is Lipschits continuous then we recover results of
condensation type when $0<\ga\ll1$.

To study the opposite case when $0<1-\ga\ll1$ one needs to apply
the transfer operators technics for suitable Banach spaces of
signed measures. These matters will be discussed in
Section~\ref{s:soft2} using a multidimensional version of the
Lasota-Yorke inequality. The idea here is to consider the action
of interactions as weak perturbations to the dynamics of each
of local systems. This perturbative approach is well known but
in our case there is an additional problem related to the fact
that in distinction to the already studied situations the transfer
operator related to the interaction might make a large contribution
to the strong norm dpending on the total number of particles in
the configuration.

One might argue that if $\ep=1$ and $\ga\ll1$ then we recover
the usual weakly intearcting CML setup. The assumption that
$\ep\ll1$ makes the already weak interactions to occur
very rare and thus should preserve the ``almost'' direct
structure of the invariant measure. Unfortunaly this ``soft''
argument does not work because the local structure of interactions
leads to the creation of small regions where the ``variation''
of a measure might become arbitrary large under the action of
interactions (see Section~\ref{s:soft2}).

Additionally to the closeness in space the interaction might
depend on the closeness in the ``lattice position''. To be precise
one considers the graph of interactions and assumes that it is
locally finite, i.e. the degree of each vertex is finite. The
degrees need not to be uniformly bounded. This allows to consider
infinite systems, discussed briefly in Section~\ref{s:gen} (see
also a mean field approach in this Section).

A few words about notation. We use a convention that for a Borel
set $A$ and a (signed/complex) measure $\mu$ the restriction
of the measure $\mu$ to the set $A$ is denoted by $\mu_{|A}$ and
$\mu(\phi):=\int\phi~d\mu$. The bar notation $\vx$ is used to mark
variables describing the CML. Note also that $|\cdot|$ is used
for very different objects throughout the paper and we follow
the convention that in the case of a subset of integers $|J|$
means its cardinality, in the case of an interval - its length,
in the case of a function $|\phi|:={\rm ess}\sup_x|\phi(x)|$, and in the
case of a (signed/complex) measure $\nw{\mu}$ stands for its weak
norm, while $\ns{\mu}$ stands for the strong norm
(see Section~\ref{s:caricature}).

\section{Direct product systems: basic ergodic constructions}\label{s:defs}

Here we give a short description of standard definitions an
constructions from ergodic theory which are necessary for the proof
of our results.

Recall that a measure $\mu$ is $\map$-{\em invariant} if and only
if $\mu(\phi\circ\map)=\mu(\phi)$ for any $\mu$-integrable
function $\phi:X\to\IR^1$.

A measurable function $\phi:X\to\IR^1$ is called {\em invariant}
with respect to a dynamical system $(\map,X,\cB,\mu)$ (or simply
$\map$-invariant), if $\phi=\phi\circ\map$ almost everywhere with
respect to the measure $\mu$.

A dynamical system $(\map,X,\cB,\mu)$ is {\em ergodic} if each
$\map$-invariant function is a constant $\mu$-a.e.

A dynamical system $(\map,X,\cB,\mu)$ is {\em weak mixing} if %
$$ \frac1n\sum_{k=0}^{n-1}|\mu(\map^{-k}A\cap B) - \mu(A)\mu(B)|
   \toas{n\to\infty}0 \qquad \forall A,B\in\cB .$$

A {\em direct product} of a pair of dynamical systems
$(\map',X',\cB',\mu')$ and $(\map'',X'',\cB'',\mu'')$ is a new
dynamical system %
$(\map'\otimes\map'',X'\otimes{X''},\cB'\otimes\cB',\mu'\otimes\mu'')$,
where the map $\map'\otimes\map'':X'\otimes{X''}\to X'\otimes{X''}$
is defined by the relation
$\map'\otimes\map''(x',x''):=(\map'x',\map''x''),$  
while all other objects are standard direct product of spaces,
$\sigma$-algebras and measures respectively.

By $A^N$ we denote the direct product of $N$ identical sets
$A\in\cB$, and by $(\map^{\otimes{N}},X^N,\cB^N,\mu^N)$ -- the
direct product of $N$ identical copies of a dynamical system
$(\map,X,\cB,\mu)$.

\begin{theorem}\label{t:mixing} Let a dynamical system
$(\map,X,\cB,\mu)$ satisfy the weak mixing property. Then for any
positive integer $N<\infty$, measurable set $\v{A}\in\cB^N$ with
$\mu^N(\v{A})>0$ and for almost any (with respect to the measure
$\mu^N$) collection $\vx:=\{x_1,x_2,\dots,x_N\}\in X^N$ there exists
a moment of time $t\ge0$ such that $\vmap^t\vx\in\v{A}$.
\end{theorem}
\proof Let a dynamical system $(\tau,Y,\cB_Y,\nu)$ be ergodic.
Then for any pair of measurable sets $A,B\in\cB_Y$ with
$\nu(A)\nu(B)>0$ there exists a positive integer
$\kappa=\kappa(A,B)<\infty$ such that $\tau^\kappa A\cap
B\ne\emptyset$. Indeed, assume that this is not true, i.e.
$\tau^nA\cap B=\emptyset$ for any positive integer $n$. Consider a
measurable set %
$A_\infty:=\bigcup_{n\in\IZ_+}\tau^nA.$ 
Then $\nu(A_\infty)\ge\nu(A)>0$ and $A_\infty\cap B=\emptyset$.
Therefore the indicator function of a measurable set $A_\infty$ of
positive $\nu$-measure is $\tau$-invariant but is not a constant
a.e. which contradicts to the ergodicity.

Therefore it is enough to show that the dynamical system
$(\map^{\otimes N},X^N,\cB^N,\mu^N)$ is ergodic. For that we shall
take advantage of the fact that the weak mixing property is
preserved under the direct product of weak mixing dynamical systems
(see e.g. \cite{Sinai}). To complete the proof its remains to note
that the weak mixing implies ergodicity. \qed%

It is of interest that one cannot weaken the conditions of
Theorem~\ref{t:mixing} replacing the weak mixing of the original
dynamical system by the ergodicity. The problem is that the direct
product of ergodic dynamical systems needs not to be ergodic as
well: consider a direct product of two identical irrational unit
circle rotations. On the other hand, the weak mixing condition is
not necessary as well (see \cite{B-raw} for an example of a
nonergodic dynamical system for which every open neighborhood
$\v{A}$ of the diagonal of $X\times X$ satisfies the claim of
Theorem~\ref{t:mixing}).

\bigskip

In what follows we often deal with the $\ep$-neighborhood of the
diagonal $\v{D}$ of $\vX\equiv X^N$ which we denote by $\v{D}_\ep$. It is
not difficult to calculate the Lebesgue measure of this set but
the product measure $\vmu(\v{D}_\ep)$ may vary sensitively with $\mu$.
To estimate it we use the techniques of measurable partitions.

Recall that a measurable partition of $(X,\cB)$ is a collection
$\Delta:=\{\Delta_i\},~\Delta_i\in\cB$ such that
$\Delta_i\cap\Delta_j=\emptyset,~\cap_i\Delta_i=X$. The diameter
of a partition is the largest diameter of its elements.

\begin{lemma}\label{l:diag} For each $0<\ep$ and a probabilistic
measure $\mu$ we have $0<\vmu(\v{D}_\ep)\le1$, moreover
$\sup_\mu\vmu(\v{D}_\ep)=1$. Assume now that $\forall \ep>0$ there
exists an partition of $X$
of diameter $\ep$ with cardinality $n_\ep\le C\ep^{-d}$. Then %
$\inf_\mu\vmu^N(\v{D}_\ep)\ge(C\ep^{-d})^{-N+1}$.
\end{lemma}
\proof Let $\Delta$ be a finite partition of $X$ of diameter
$\ep>0$. Since $X$ is compact such partitions always exist. Denote
$a_i:=\mu(\Delta_i)$ then $\sum_ia_i=\mu(X)=1$.
Observe now that %
$\Delta_i^N:=\underbrace{\Delta_i\times\dots\times\Delta_i}_N
          \subset\v{D}_\ep$. %
Therefore the lower estimate of $\vmu(\v{D}_\ep)$ follows from the
trivial inequality
$\vmu(\v{D}_\ep)\ge\sum_i(\mu(\Delta_i))^N=\sum_ia_i^N>0$.

Now let $\mu$ be concentrated at a single point $x\in X$, i.e.
$\mu(\{x\})=1$. Then obviously
$\vmu(\v{D}_\ep)\equiv1~\forall\ep>0$.

It remains to prove the lower estimate under the additional
assumption about the cardinality of the partition $\Delta$.
Consider a power average %
$S_N(\{a_i\}):=\left(\frac1{n_\ep}\sum_{i=1}^{n_\ep}a_i^N\right)^{1/N}$
of $n_\ep$ nonnegative entries $a_i$. It is known that
$S_N(\cdot)$ is monotonous on the nonnegative parameter $N$. Thus
$S_N(\{a_i\})\ge S_1(\{a_i\})\equiv\frac1{n_\ep}$. Therefore %
$\vmu(\v{D}_\ep)\ge\sum_i(\mu(\Delta_i))^N=\sum_ia_i^N
 =n_\ep (S_N(\{a_i\}))^N\ge n_\ep(S_1(\{a_i\}))^N =n_\ep^{-N+1}$. \qed

\begin{lemma}\label{l:nonsingular} Let $\v{A}\subset\vX$ with
$\vmu(\v{A})=0$. Then
$\vmu((\vmap)^{-1}\v{A})+\vmu((\vQ_\ep)^{-1}\v{A})=0~~\forall
\ep,\gamma\ge0$.
\end{lemma}%
\proof The fact that $\vmu((\vmap)^{-1}\v{A})=0$ follows from the
nonsingularity of the local map $\map$. If $\ep=0$ the interaction
does not occur and hence $\vQ_0$ is an identical map. Therefore it
is enough to consider $\ep>0$. $\forall N\in\IZ_+$ and $\forall
\vx\in\vX^N$ the interaction map $\vQ_\ep$ may be written as a
nonsingular linear map. However the matrix defining this map
depends on the ``grouping'' $J_i(\vx)$. Nevertheless for the
finite number $N$ of particles the total number of various
``grouping'' is finite and hence $\vQ_\ep$ may have only a finite
number of nonsingular representations. \qed

A map $\map:X\to X$ induces the {\em transfer operator} $\map^*$
acting in the space of signed measures (generalized functions)
$\cM$ on $X$ by the formula
$\map^*\mu(A):=\mu(\map^{-1}A)$ for each $A\in\cB$ and
$\mu\in\cM$. From this point of view a measure $\mu$ is
$\map$-invariant if and only if $\map^*\mu=\mu$.

A probabilistic measure $\mu_\map\in\cM$ is called the
{\em Sinai-Bowen-Ruelle} (SBR) measure for the dynamical
system $(\map,X,\cB)$ and a reference measure $m$
(say Lebesgue measure on $X$) if there exists an
open subset $Y\subseteq X$ such that for each probabilistic
measure $\mu\in\cM$ absolutely continuous with respect to $m$
and such that $\mu(Y)=1$ we have weak convergence %
$\frac1n\sum_{t=0}^{n-1}{\map^*}^t\mu\toas{n\to\infty}\mu_\map$.
The set $Y$ is called the {\em basin of attraction} for the
measure $\mu_\map$. Obviously an SBR measure is $\map$-invariant.

There are also different approaches for the definition of the SBR
measure and we refer the reader to \cite{BB} for their discussion
and conditions under which those approaches agree with each other.

\section{Rigid interactions}\label{s:rigid}

\subsection{Proof of Theorem~\ref{t:sync2}}

Let $N=2$. Consider the $\ep$-neighborhood of the diagonal
$\v{D}_\ep$ in $\vX$. By Lemma~\ref{l:diag} $\vmu(\v{D}_\ep)>0$.
Therefore by Theorem~\ref{t:mixing} for $\vmu$-a.a. $\v{x}\in\vX$
there exists the first moment of time $0\le t(\v{x})<\infty$ such
that $\vmap^{t(\v{x})}\v{x}\in\v{D}_\ep$. Denote the set of full
$\vmu$-measure for which this holds by $\vX_\ep$.

By the definition of the rigid interaction,
$\vQ_\ep\v{D}_\ep\subseteq\v{D}$ while $\vQ_\ep\v{x}\equiv\v{x}$
if $\v{x}\notin\v{D}_\ep$. Therefore $\forall \v{x}\in\vX_\ep,
~~t\in\{0,\dots,t(\v{x})-1\}$ we have
$\vmap_\ep^t\v{x}\equiv\vmap^t\v{x}$ and
$\vmap_\ep^{t(\v{x})}\v{x}\in\v{D}$ which proves the global
synchronization if $N=2$.

Local synchronization for an arbitrary $N\ge2$ follows from the
invariance of the $\ep$-neighborhood of the diagonal with respect to
the dynamics. However, if $N\ge3$ the observation that
$\vmap_\ep^t\v{x}\equiv\vmap^t\v{x}~~\forall
t\in\{0,\dots,t(\v{x})-1\}$ does not hold.

It remains to show that when $N\ge3$ even the analyticity of the map
$\map$ does not guarantee the global synchronization.

\begin{lemma}\label{l:nosync} Let $X:=S^1$ (unit circle), the
local system be governed by the doubling map $\map x:=\{2x\}$ and
let $N=3$. Then $\forall~ 0<\ep\ll1$ the desynchronization with the
domain of positive Lebesgue measure occurs.
\end{lemma}%
\proof Let $x_i\in X, ~i\in\{1,2,3\}$ and denote $a:=x_2-x_1,~
b:=x_3-x_2$. Denote by $A$ the subset of $\vX$ for which
$0<a,b\le\ep$ and $a+b>\ep$. Under the interaction the coordinates
$x_i$ will be changed to $x_i'$ such that the distances between them
will become equal to $a':=(a+2b)/6$ and $b':=(2a+b)/6$. Applying the
doubling map we get the new pair of distances $a'':=(a+2b)/3$ and
$b'':=(2a+b)/3$. Since $0<a'',b''\le\ep$ and $a''+b''=a+b>\ep$ the
new configuration again belongs to the set $A$. Now the observation
that the product Lebesgue measure $\v{m}(A)=\ep^2/2>0$ finishes the
proof. \qed

The extension of Lemma~\ref{l:nosync} for the case when dim$(X)>1$
is straightforward. Moreover, the discussion after the formulation
of Theorem~\ref{t:sync2} demonstrates that for large $N$ the local
map demonstrating the desynchronization may be chosen
$o(1/N)$-close to identical for each $\ep>0$.

\subsection{Proof of Theorem~\ref{t:sbr2}}

One might think that Theorem~\ref{t:sbr2} is a direct consequence
of Theorem~\ref{t:sync2}. Indeed by Theorem~\ref{t:sync2} a.a.
trajectories of our CML $\forall\ep>0$ after a finite number of
iterations hit the diagonal. Hence any forward invariant set of
$\vmap_\ep$ belongs to the diagonal $\v{D}$ and obviously the
analysis of the invariant measures may be restricted to the
``one-dimensional'' dynamics on the diagonal. On the other hand,
we do not assume that the $\map$-invariant measure $\mu$ is unique
and hence typically there is a subset $Y\subset X$ on untypical
points leading to statistics different from $\mu$. Still the
$\mu$-measure of this exceptional set is zero. Therefore using
that number of iterations before to hit the diagonal is finite
(but not uniformly bounded) and the result of
Lemma~\ref{l:nonsingular} we deduce that the set of
$\vmu$-``typical'' points is of full $\vmu$-measure.

\subsection{The closest rigid interactions}\label{s:closest}

As usual only the particles from the $\ep$-neighborhood of $x_i$
will be included to $J_i(\vx)$, but now we consider a special
(simplified) choice of %
$$J_i(\vx):=\{j:~~ \rho(x_i,x_j)\le\ep, ~~
              \rho(x_i,x_j)=\min_{x_k\ne x_i}\rho(x_i,x_k)\} ,$$
i.e. this collection contains only the {\em closest} particles to
the $i$-th one. We shall refer to this choice of $J_i(\vx)$ as the
{\em closest interaction}. Note that the set of configurations
$\vx$ for which $\max_i|J_i(\vx)|>2$ has Lebesgue measure 0.
\begin{theorem}\label{t:syncN-simple} Let the $(\map,X,\cB,\mu)$ be a
weakly mixing measurable DS and $\mu$ be its only SBR measure.
Assume also that the interaction is rigid with the above choice of
$J_i(\vx)$. Then for each $N\in\IZ_+, ~\ep>0$ the $\v\mu$-global
synchronization takes place.
\end{theorem}%
\proof For a configuration $\vx\in\vX$ denote by $\ell(\vx)$ the
minimum distance between particles in $\vx$. Observe that for
$\v\mu$-a.a. configurations each $J_i(\vx)$ consists of at most
one particle and the minimal distance is achieved at a single pair
of particles. By the definition of the interaction if
$\ell(\vx)\le\ep$ these two particles may interact only between
themselves and hence their positions after the interaction will
coincide with the common center of gravity.

On the other hand, if $\ell(\vx)>\ep$ no interactions occur and
one may use the same argument as in the proof of
Theorem~\ref{t:sync2} to show that $\v\mu$-a.s. this effect takes
place.

To finalize the inductive construction, observe that for each
particle the matching with some other particles takes place after
a finite number of time steps and hence by
Lemma~\ref{l:nonsingular} the $\v\mu$ measure of ``non generic''
initial configurations leading to the non-uniqueness
of $\ell(\vx)$ is zero. \qed%

\section{Soft interactions. Case $0\le\ga\ll1$}\label{s:soft1}

In this Section we study the intermediate case when $0<\ga<1$.
It is easy to see that if $\ga$ is close enough to $1$ then we
are basically in the same situation as in the case of rigid
interactions.

\begin{theorem}\label{t:soft-rigid} Let the $(\map,X,\cB,\mu)$ be a
weakly mixing measurable DS and $\mu$ be its only SBR measure.
Assume also that there is a constant $0<\La<1/\ga$ such that
$\rho(\map{x},\map{y})\le\La\rho(x,y)~~\forall x,y\in X$
Then the claims of Theorems~\ref{t:sbr2} and
\ref{t:sync2} hold true.
\end{theorem}
\proof Let $\vx\in\v{D}_\ep$. Then $\rho(x_i,x_j)\le\ep~~\forall i,j$ and
hence all ``particles'' do interact with each other. Denote by $z$ their
common center of gravity. Then $(\vQ_\ep\vx)_i:=\ga x_i+(1-\ga)z$ which
implies that
$$\rho((\vQ_\ep\vx)_i , (\vQ_\ep\vx)_j)
  \le \ga\rho(x_i,x_j)~~\forall i,j .$$
Thus $\vQ_\ep\v{D}_\ep\subseteq \v{D}_{\ga\ep}$.

On the other hand $\vmap \v{D}_\ep \subseteq \v{D}_{\La\ep}$
by the assumption on the map $\map$. Therefore %
$$ \vmap\circ\vQ_\ep \v{D}_\ep \subseteq \v{D}_{\La\ga\ep}
  \subseteq \v{D}_\ep .$$
Moreover, %
$$ \rho((\vmap_\ep^t\vx)_i,(\vmap_\ep^t\vx)_j)
 \le(\La\ga)^t\rho(x_i,x_j)
 \toas{t\to\infty}0~~\forall i,j .$$

The completion of the proof follows exactly to the same arguments
as in the case of rigid interactions. \qed %

\begin{corollary} Let $\La_\map$ be the modulus of the largest
Lyapunov multiplier of the map $\map$. Then the claims of
Theorems~\ref{t:sbr2} and \ref{t:sync2} hold true if and only if
$\ga<\ga_0:=1/\La_\map$.
\end{corollary}%
\proof The direct statement follows from the argument above
applied to $\v{D}_{\ep/\La_\map}$ rather than to $\v{D}_\ep$. To
prove the inverse statement one observes that the modulus of the
largest Lyapunov multiplier of the map $\vmap\circ\v{Q}_\ep$
cannot be smaller than $\La_\map\ga>1$. \qed%

\n{\bf Remark.} Despite the similarity between the case under consideration
and the case of rigid interactions there is an important difference in
that once $\vx\in\v{D}_\ep$ all ``particles'' will immediately form
a cluster on the next time step in the rigid case, while an infinite number
of iterations is nessasary for this in the soft case.

\section{Soft interactions. Case $0<1-\ga\ll1$}\label{s:soft2}

Since $\ga=1$ correspond to the absence of interactions (and thus the
SBR measure of the multicomponent system is equal to the direct product
of local SBR measures) one expects to observe a kind of phase transition
when the parameter $\ga$ grows from 0 to 1. In what follows we shall
study what happens when the parameter $\ga$ becomes very close to 1.


In the Introduction we already mentioned that for our purpose it is
enough to assume that the local transfer operator $\map^*$ be
quasi-compact in a suitable Banach space of signed measures. This
means that $\map^*$ may be represented as a sum of a contraction
and a compact operator. Below we shall show that this property
is satisfied e.g. for the so called piecewise-expanding maps.

To this end we need to introduce a proper Banach space of signed
measures and to describe their properties.

\subsection{Transfer-operator approach and BV measures}\label{s:caricature}

Recall that $\vX$ is a unit $N$-dimensional Eucledian cube equipped
with the standard Borel $\sigma$-algebra $\v\cB\equiv\cB^N$.
The map $\vmap:\vX\to\vX$ induces the {\em transfer operator}
$\vmap^*$ acting in the space of signed measures (generalized functions) $\v\cM$ on $\vX$ by the formula
$\vmap^*\mu(\v{A}):=\mu(\vmap^{-1}\v{A})$ for each $\v{A}\in\v\cB$ and
$\v\mu\in\v\cM$.

Starting from \cite{BKL,Bl3} the approach to the analysis of
transfer operators in terms of the so called dual norms proved to
be efficient and became popular. To introduce the dual norms in
the space of signed measures (generalized functions) on $\vX$ we
start with spaces of test-functions %
\bea{ \cF\a:=\{\phi\in C^1(\vX):~~|\phi|\le1\}, \quad
   \cF_0:=\{\phi\in C^1(\vX):~~|\phi|\le1,~ \phi|_{_{\partial \vX}}=0\}
   ,\\
   \cF_L\a:=\{\phi\in \IL(\vX):~~|\phi|_\infty\le1\} .} %

The following are two versions of the ``variation'' of a signed
measure:
$$ \V(\vmu):=\max_i\sup_{\phi\in\cF}\vmu(\partial_i\phi), \qquad
   \V_0(\vmu):=\max_i\sup_{\phi\in\cF_0}\vmu(\partial_i\phi) .$$
The latter functional gives the variation of the density of the
measure $\mu$ with respect to the Lebesgue measure, while the
former also gives the variation of this density but considered as
a function from $\IR^N$ taking zero value outside of $X$. An
important advantage of the above definition of the variation is
that in the case of a measure having the direct product structure
properties of its variation can be easily obtained from their
one-dimensional counterparts. Therefore we shall give proofs only
for one-dimensional statements and refer the reader e.g. to
\cite{KL05} for the multidimensional setting.

Define also the $L^1$-norm of the (signed) measure:
$$|\vmu|:=\sup_{\phi\in\cF_L}\vmu(\phi)$$ which we shall call
the {\em weak} norm.

\begin{lemma}\label{l:var-properties}
(a)  $\V(\vmu_{|\v{Y}})\le\V(\vmu)$; \par
(b) $\V_0(\vmu_{|\v{Y}})\le\V(\vmu_{|\v{Y}})
    \le2\V_0(\vmu_{|\v{Y}})+2|\vmu_{|Y}|/{\vm(\v{Y})}$; \par
(c) if $\v{Y}$ is a proper rectangle\footnote{i.e. a
       direct product of $n$ intervals.} then %
    $|\vmu_{|Y}|\le\frac12~\vm(\v{Y})~\V(\vmu)$,
    in particular, $|\vmu|\le\frac12\V(\vmu)$.
\end{lemma}%
 \proof (a) $\V(\mu_{|Y})=\sup_{\phi\in\cF}\mu_{|Y}(\phi)
               \le\sup_{\phi\in\cF}\mu(\phi)=\V(\mu)$.

(b) For $\phi\in\cF$ set
    $\phi_0:=x(\phi(x)-\phi(1))+(1-x)(\phi(x)-\phi(0))$.
Then $|\phi_0|\le2|\phi|$ and since $\phi_0(0)=\phi_0(1)=0$ we
obtain that $\frac12\phi_0\in\cF_0$. Therefore
$\mu(\phi')=\mu(\phi_0+(\phi(1)-\phi(0)))\le2\V_0(\mu)+2|\mu|$,
which proves the inequality for the case $Y=X$. The general case
can be proven similarly.

(c) Decompose the signed measure $\mu_{|Y}:=\mu_+-\mu_-$ into positive
and negative components and set $Y_\pm:={\rm supp}(\mu_\pm)$. The function
$\phi(x):=m(Y_+\cap[0,x]) - m(Y_-\cap[0,x]) - \frac12m(Y)$ is continuous
on $X$ and $|\phi|\le\frac12 m(Y)$. On the other hand, by definition
$|\mu_{|Y}|=\mu(\phi')\le\frac12 m(Y) \V(\mu)$ since $2\phi/m(Y)\in\cF$
is a valid test-function.
\qed   %

Therefore the functional $\V(\vmu)$ is actually a norm (which we denote
by $\ns{\vmu}$) and is equivalent to a more common strong norm
$\V_0(\vmu)+\nw{\vmu}$. Therefore we shall refer to $\ns{\cdot}$ as a
{\em strong norm}. Note also that for the Lebesgue measure $\V(\vm)=2$.
The set of (signed) measures $\mu$ with $\ns{\vmu}<\infty$ we shall call
measures of {\em bounded variation} and denote this set by ${\bf BV}$.

Using this terminology we may rewrite Theorem~\ref{t:main2} claiming
the convergence to the direct product measure as follows.

\begin{theorem}\label{t:stable} Let the map $\map$ have the only one
SBR measure $\mu_\map$, and let there are constants
$0\le\theta<1\le\Theta<\infty$ such that %
\beq{e:ly-a}{ \ns{\map^*\mu}\le\theta\ns{\mu}+\Theta\nw{\mu} } %
for each $\mu\in${\bf BV}. If $2^{N}\theta/\ga<1$
then for each $0\le\ep\ll1$ the CML $(\vmap_\ep,\vX)$ has the only
one SBR measure $\vmu_\ep\toas{\ep\to0}\vmu_\map$ -- the direct
product of the local SBR measures.
\end{theorem}

To give a specific model satisfying to our assumptions consider the
class of piecewise expanding maps. Let $X:=[0,1]$ and
$\tau:[0,1]\to[0,1]$ be a piecewise $C^2$-smooth map, i.e. there
is a finite partition of $X$ into intervals $X_i$ on each of which
the map $\tau$ is bijective, $C^2$-smooth, and
$\inf_{x}|\tau'(x)|\ge\la>0$, ~$\beta_1(\tau):=\frac2{\la \min_i|X_i|}$,
~$\beta_2(\tau):=\sup_x|(1/\tau'(x))'|<\infty$. Such maps are
called $\la$-{\em expanding}. Set
$\beta(\tau):=\beta_1(\tau)+\beta_2(\tau)$.

\begin{lemma}\label{l:LY} (Lasota-Yorke inequality) Let the maps
$\tau_1,\tau_2,\dots,\tau_N$ be $\la$-expanding, $\v\tau$
stands for their direct product, and let
$\beta(\v\tau):=\max_i\beta(\tau_i)$. Then %
\beq{e:LY}{\ns{\v\tau^*\vmu}
           \le\frac2\la\ns{\vmu} + \beta(\v\tau)~\nw{\vmu} .}%
\end{lemma}
\proof First observe that from Lemma~\ref{l:var-properties},(b) it
follows that for each $\phi\in\cF$ %
\beq{e:decomp}{\mu(\phi') \le 2|\phi|
         \left(\V_0(\mu) + \frac1{|X_i|}~|\mu|\right) }%
As $(\phi\circ\tau)'(x)=\phi'(\tau(x))\cdot\tau'(x)$ for each
$x\in X\setminus(\cup_i\partial X_i)$, we have %
$$ \tau^*\mu(\phi') = \mu(\phi'\circ\tau)
 = \mu((\phi\circ\tau)'/\tau')
 = \mu(((\phi\circ\tau)/\tau)')
   - \mu((\phi\circ\tau)\cdot(1/\tau)') .$$
To estimate the first term we apply (\ref{e:decomp}), while the
second term is bounded by $\beta_2(\tau)~|\mu|$. \qed %

\begin{corollary}\label{c:inv-mes} Under the assumptions of
Lemma~\ref{l:LY} there exists a probabilistic $\v\tau$-invariant
measure $\mu_{\v\tau}$.
\end{corollary}
\proof Choose $k\in\IZ_+$ large enough such that $\la^k>2$.
Denote $r:=2/\la^k<1$ and let $\vmu\in{\bf BV}$ be a probabilistic
measure. Then for each $n\in\IZ_{+}$ we have: %
$$\ns{{\v\tau^*}^{nk}(\vmu)} \le r^{n}\ns{\vmu}
                          + \frac{\beta(\v\tau)}{1-r}\nw{\vmu} .$$
Thus the sequence $\mu_{n}:={\v\tau^*}^{nk}\vmu$ satisfies the
conditions of the embedding of {\bf BV} to the set of measures
having absolutely continuous densities with respect to Lebesgue
measure, which we denote by $\IL^1$, and, hence, there exists
the limit point of this sequence, i.e. %
$\vmu_{n_{i}}\toas{i\to\infty}\vmu_{\infty}\in{\bf BV}$ with
$\ns{\mu_{\infty}} \le \frac{\beta(\v\tau)}{1-r}\nw{\mu}$.
On the other hand, being a limit point of this sequence the
measure $\v\mu_{\infty}$ satisfies the relation
$\v\tau^*\vmu_{\infty}=\vmu_{\infty}$ and thus
is $\v\tau$-invariant. \qed%

\begin{corollary}\label{c:pe} The transfer operator corresponding
to the map $\v\tau$ under an additional assumption of the uniqueness
of SBR measures for each map $\tau_i$ satisfies the conditions of
Theorem~\ref{t:stable} with $\theta:=2/\la$ and $\Theta:=\beta(\v\tau)$.
\end{corollary}

\Bfig(320,140)
      {\footnotesize{
       \put(0,0){\bpic(150,140){
       \put(0,0){\vector(1,0){150}} \put(0,0){\vector(0,1){150}}
       \bezier{100}(75,0)(75,75)(75,150) \bezier{100}(0,75)(75,75)(150,75)
       \bezier{100}(60,0)(60,75)(60,150) \bezier{100}(90,0)(90,75)(90,150)
       \thicklines
       \bline(10,10)(1,1)(50) \bline(150,150)(-1,-1)(60)
       \bline(75,75)(3,1)(15) \bline(75,75)(-3,-1)(15)
       \put(61,80){$\al_i x + c_i$}
       \put(145,-10){$x$}  \put(-12,140){$\tau x$}
       \put(60,-10){$a_i$} \put(90,-10){$a'_i$}
       }}
       \put(180,0){\bpic(150,140){
       \put(0,0){\vector(1,0){160}} \put(0,0){\vector(0,1){150}}
       \bezier{100}(60,0)(60,75)(60,130) \bezier{100}(90,0)(90,75)(90,130)
       \bezier{200}(0,30)(60,180)(145,20)       
       \bezier{100}(60,120.5)(77,125)(90,110) 
       \thicklines
       \bline(60,120.5)(-1,0)(60) \bline(90,110)(1,0)(55)   
       \put(60,-10){$a_i$} \put(90,-10){$a'_i$} \put(155,-10){$x$}
       \put(75,123){$\Phi_i$} \put(75,91){$\phi$}
       }}
      }}
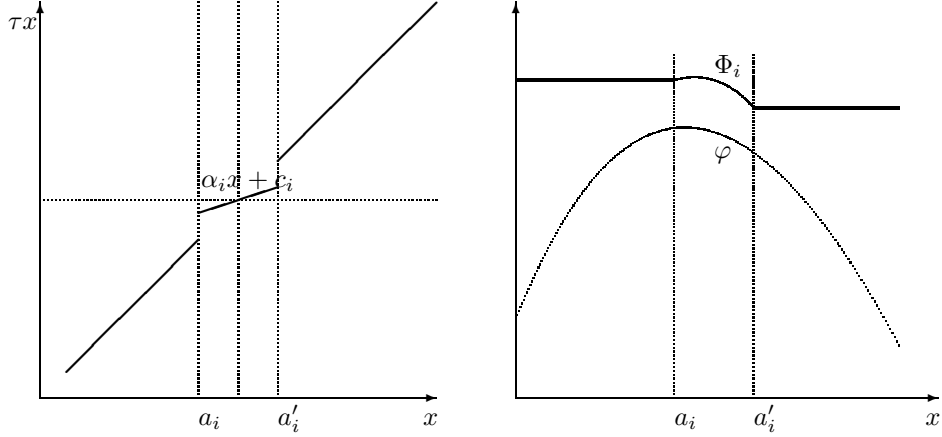
{Local structure of the map $\tau$ (left). ~ Decomposition of the test-function (right).
\label{f:graph}}

Consider now a more specific family of piecewise $C^2$-smooth maps
$\tau:X\to X$. For a given positive integer $n$ let
$0\le a_1<a_1'\le a_2<a_2'\le\dots\le a_n<a_n'\le1$
and set $A_i:=[a_i,a'_i], ~i\in\{1,2,\dots,n\}$, $A:=\cup_iA_i$,
and $B:=X\setminus A$. These intervals define a partition of $X$.
Then we define %
$\tau x:=\function{\tau_i x:=\al_i x+c_i &\mbox{if } x\in A_i \\
                    x                   &\mbox{otherwise}}$,
where $\tau_iA_i\subseteq A_i, ~\al_i>0~~\forall i$.
Fig.~\ref{f:graph} (left) demonstrates the shape of $\tau$ in a
neighborhood of an interval $A_i$.

We shall be interested in the properties of the transfer-operator $\tau^*$
when $|A|\ll1$. Due to this restriction the application of
Theorem~\ref{l:LY} gives an estimate with the second term going
to infinity as the diameter of the partition goes to zero. In order
to overcome this difficulty we develop a new approach to estimate
the norm of the transfer operator in this case.

\?{Let $X:=[0,1]$ and $\tau:[0,1]\to[0,1]$ be the map depicted in
Fig.~\ref{f:graph} (left) and defined by the relation %
$\tau x:=\function{\al x+c &\mbox{if } x\in A \\
                    x      &\mbox{otherwise}}$, %
where $a,b\in(0,1), A:=[a,b], |A|\le\ep, 2c=(1-\al)(a+b)/2$. }

\begin{lemma}\label{l:discontinuous-pert}
$\ns{\tau^*\mu}\le(n+1+\sum_i\frac1{\al_i})~\ns{\mu}$, and %
$\nw{\mu-\tau^*\mu}\le \frac12 (n+2+\sum_i\frac1{\al_i})~m(A)~\ns{\mu}$.
\end{lemma} %
\proof The idea used in the proof of
Lemma~\ref{l:var-properties},(b) is to interpolate linearly
between the values of the test-function at boundary points of the
partition $\{X_i\}$ and to estimate the contribution of this
interpolation into the integral against the weak norm, rather than
the strong one. In the case under consideration the lengths of the
intervals of monotonicity might be arbitrary small which does not
allow to apply his trick. Instead we shall treat each interval of
monotonicity separately extending the test-function by two
constants equal to the values at boundary points outside of the
interval (see Fig.~\ref{f:graph} (right)).

Observe that for each $i$ the function $\tau_i(x)$ can be extended
as a linear function to the whole $X$. For a test-function
$\phi\in\cF$
and a (signed) measure $\mu$ we have %
\bea{  \tau^*\mu(\phi') \a= \mu(\phi'\circ\tau)
 = \mu_{|B}(\phi') + \sum_i\mu_{|A_i}(\phi'\circ\tau) \\
 \a= \mu_{|B}(\phi') + \mu_{|A}(\phi') - \mu_{|A}(\phi')
   + \sum_i\mu_{|A_i}(\phi'\circ\tau) \\ %
 \a= \mu(\phi') + \sum_i^n\left(\mu_{|A_i}(\phi'\circ\tau)
                            - \mu_{|A_i}(\phi') \right) .} %
Let $x\in A_i$ then
$\phi'\circ\tau=(\phi\circ\tau)'\cdot(\tau_i')^{-1}
=\frac1{\al_i}(\phi\circ\tau)'$ and $\phi\circ\tau\in\cF$ is a valid
test-function. Therefore %
$\mu_{|A_i}(\phi'\circ\tau)\le\frac1{\al_i}\V(\mu)$. Summing up all
contributions from the integrations over $A_i$ we get
$$\tau^*\mu(\phi') \le \left(1+n+\sum_i\frac1\al_i\right)~\ns{\mu} $$
since $\mu_{|A_i}(\phi')\le\ns{\mu}~~\forall i$ by
Lemma~\ref{l:var-properties}.

\?{In the case of $\Phi_0$ we proceed a bit differently:
$\mu(\phi'\circ\tau)
   = \mu(\phi',[0,a]\cup[b,1])
     + \mu(\Phi') - \mu(\Phi',[0,a]\cup[b,1]) 
  = \mu(\phi',[0,a]\cup[b,1])
     + \mu(\Phi') - \frac{\mu([1,a])}{\tau a}
                  - \frac{\mu([b,1])}{1- \tau b} 
  \le \V_0(\mu) + \frac1\al\V_0(\mu)
                + (\frac1{\tau a}+\frac1{1- \tau b})|\mu|.$}

To estimate $\nw{\mu-\tau^*\mu}$ observe that the measures differ
only on the intervals $A_i$. Set $\nu:=\mu-\tau^*\mu$ and
$A:=\cup_iA_i$. Then %
\bea{\nw{\nu}\a=|\nu_{|A}|\le\frac12 m(A)\V(\nu)
   \le \frac12 m(A) (\ns{\mu} + \ns{\tau^*\mu}) \\
   \a\le \frac12 m(A)\left(n+2+\sum_i\frac1{\al_i}\right)\ns{\mu}
   .}
\qed %

It might seem that the multiplier $(n+1+\sum_i\frac1{\al_i})$ is overpessimistic,
but the trivial example of the Lebesgue measure $m$ on $X$ immediately
shows that $\ns{\tau^*m}=2(n+1+\sum_i\frac1\al_i)$ while $\ns{m}=2$.

As we shall see the argument used in this proof is the key point
in the proof of Theorem~\ref{t:stable}. Apart from this
Lemma~\ref{l:discontinuous-pert} allows to apply the operator
approach to a new class of small but discontinuous perturbations.

\begin{theorem}\label{t:disc1d}
Let $\map:X\to X$ be a $\la$-expanding map with $\la>1$ having a
unique SBR measure $\mu_\map$ and let $\tau$ be a piecewise linear
map described above with $|A|=\delta$ and a given collection of
slopes $\{\al_i\}$ such that $2(n+1+\sum_i\frac1{\al_i})<\la$.
Then for each $0<\delta\ll1$ the dynamical system
$(\tau\circ\map,X)$ has a unique smooth invariant measure
$\mu_\delta\toas{\delta\to0}\mu_\map$.
\end{theorem}%
\proof Combining results of Lemmas~\ref{l:LY},
\ref{l:discontinuous-pert}
we get %
\bea{ \ns{(\tau\circ\map)^*\mu}
    \a\le (n+1+\sum_i\frac1{\al_i})~\ns{\map\mu} \\
    \a\le
    \left(n+1+\sum_i\frac1{\al_i}\right)\frac2\la\ns{\mu}
  + \left(n+1+\sum_i\frac1{\al_i}\right)\beta(\tau)\nw{\mu} .}%
Therefore we are in a position to apply Corollary~\ref{c:inv-mes}
to the map $\tau\circ\map$ which yields the existence of the
probabilistic invariant measure $\mu_\delta\in{\bf BV}$.

On the other hand, by Lemma~\ref{l:discontinuous-pert}
$$ \nw{\mu - \tau^*\mu} \le C\delta\ns{\mu} $$
which by the now standard perturbation argument (see e.g.
\cite{Bl-mon,BKL}) implies the
convergence $\mu_\delta\toas{\delta\to0}\mu_\map$. \qed %

To finish this preparatory part let us formulate an estimate of
the action of the transfer operator of a contracting affine map
which can be proven by a direct inspection.

\begin{lemma}\label{l:contraction} Let $\vX$ be the $N$-dimensional
unit cube and let $\v\tau (\vx):=G\vx+H$ be an affine map from $\vX$
into itself such that $\ell(G\xi)\ge\alpha\ell(\xi)$ for each
$\xi\in\IR^N$ and any norm $\ell(\cdot)$. Then
$\ns{\v\tau^*\vmu}\le1/\alpha\ns{\vmu}$.
\end{lemma}

\subsection{Proof of Theorem~\ref{t:stable}}

One is tempted to argue as follows. Assume that the parameter
$\ep$ is of order of the diameter of the set $X$. Then each
pair of particles is interacting between themselves and we are
coming to the well known mean field interaction model.
One can show that when $\ga$ goes to $1$ the only SBR measure
of the mean field model converges to the direct product measure,
i.e. the subsystems behave independently. Now the decrease of
$\ep$ leads only to the decrease of the frequency of interactions
between particles. Thus already ``almost'' independent particle
are becoming even more independent. In fact the situation is
much more complicated. The point is that the existence of a
large number of small islands in the phase space where different
combinatorial types of interactions actually take place
leads to a severe amplification of the variation of a measure
under the action of the operator $\Q^*$. Consider this in detail.

The definition of the map $\Q$ implies that configurations
$\vx\in\vX$ whose coordinates have pair distances larger
than $\ep$ are fixed points of the map $\Q$. On the remaining
part of the phase space $\vX$ consisting of a large (of order
$2^N$) number of disjoint components of small Lebesgue measure
the map $\Q$ is linear and contracting in each of them. Thus
the structure of the multi-dimensional map $\Q$ is very similar
to the one-dimensional map $\tau$ considered in
Lemma~\ref{l:discontinuous-pert}. Therefore we shall use
basically the same strategy for the proof.

For a given test-function $\phi\in\cF$ we need to evaluate
the functional $\max_i \Q^*\vmu(\partial_i\phi)$.

Fix some index $i$ and for each subset $J$ of different integers
belonging to the set $\{1,2,\dots,N\}$ and containing the
index $i$ define a set
$$ A_J:=\{\vx\in\vX:~~|x_i-x_j|\le\ep~\forall j\in J,~~
                      |x_i-x_k|>\ep~\forall k\notin J \}, $$
and let $A:=\cup_{|J|>1}A_J, ~B:=\vX\setminus A$.
Then the interaction with the $i$-th particle occurs only for
$\vx\in A$ and the sets $A,B$ define a finite partition of $X$.

We have %
\bea{ \Q^*\vmu(\partial_i\phi) \a= \vmu(\partial_i\phi\circ\Q)
 = \vmu_{|B}(\partial_i\phi\circ\Q)
 + \sum_{|J|>1}\vmu_{|A_J}(\partial_i\phi\circ\Q)  \\
 \a= \vmu_{|B}(\partial_i\phi)
 + \vmu_{|A}(\partial_i\phi) - \vmu_{|A}(\partial_i\phi)
 + \sum_{|J|>1}\vmu_{|A_J}(\partial_i\phi\circ\Q) \\
 \a= \vmu(\partial_i\phi) + \sum_{|J|>1}
     \left(\vmu_{|A_J}(\partial_i\phi\circ\Q)
   - \vmu_{|A_J}(\partial_i\phi)  \right) } %
Denote by $\phi'$ the vector of partial derivatives of $\phi$, by
$\Q'$ the matrix of partial derivatives of the map $\Q$, and by
$(q_{ij})$ the matrix inverse to the matrix $\Q'$ (i.e.
$(q_{ij})=(\Q')^{-1}$). Then %
\bea{ (\partial_i\phi)\circ\Q \a= ((\phi\circ\Q)\cdot(\Q')^{-1})_i
 = \sum_j\partial_j(\phi\circ\Q)\cdot q_{ji} \\
 \a= \sum_j\partial_j(\phi\circ\Q\cdot q_{ji})
 - \sum_j\phi\circ\Q\cdot \partial_jq_{ji} .}
Denoting $\phi_{ij}:=\phi\circ\vQ\cdot q_{ij}\in\IC^1$
we rewrite the expression for the action of the transfer-operator
on a measure restricted to $A_J$ as follows: %
$$ \vmu_{|A_J}(\partial_i\phi\circ\Q)
 \le \sum_j \sup_x|\phi_{ji}(x)|\cdot\V(\vmu)
   + \sum_j \sup_x|\partial_j q_{ji}(x)|\cdot\nw{\vmu} .$$
The interaction inside of each region $A_J$ is described
by a linear function and thus the last term is equal to zero
for $\vx\in A_J$.

In general the prefactor $\sum_j \sup_x|\phi_{ji}(x)|$ can be
estimated as %
$$ \sum_j \sup_x|\phi_{ji}(x)|
\le \sum_j \sup_x|((\Q')^{-1})_{ji}| .$$ %
In our case the upper estimate can be done explicitly. Observe
that for a given set $A_J$ the restriction of the map $\Q$ to
$A_J$ is a linear map (which we denote by $L$) has a very simple
structure. Namely, considering only ``interacting coordinates''
one can rewrite this map as an affine map $L\vy:=G\vy+H$ with
$G:=\ga+\frac{1-\ga}{n}E$. Here $n$ stands for the number of the
``interacting coordinates'' and all entries of the $n\times n$
matrix $E$ are ones. It is easy to check that
$\ell(G\xi)\ge\gamma\ell(\xi)$ for any norm $\ell(\cdot)$ and
$\xi\in\IR^n$ and the equality is achieved on a vector $\xi$
having the only one nontrivial coordinate. Therefore using
Lemma~\ref{l:contraction} and setting $\alpha:=\gamma$ we get
$||L^*\vmu||\le1/\gamma$. Thus the prefactor can be estimated from
above as $1/\gamma$ uniformly on $J$.

Now using that the number of different collections $J$
cannot exceed $2^N$ we get
$$ \ns{\Q^*\vmu} \le \frac{2^N\theta}{\ga} \ns{\vmu} .$$
Combining this result with the Lasota-Yorke type
inequality~(\ref{e:ly-a}) for the direct product map we estimate
the strong norm of the action of the transfer-operator of the CML as
$$ \ns{\vmap_\ep^*\vmu}
  \le \theta\ns{\Q^*\vmu} + \Theta\nw{\Q^*\vmu}
  \le \frac{2^N\theta}{\ga}\ns{\vmu} + \Theta\nw{\vmu} .$$
Now using again the same trick as in the proof of
Lemma~\ref{l:discontinuous-pert} we show that the measures
$\vmap_\ep^*\vmu$ and $\vmap^*\vmu$ are close in the weak
($\IL^1$) norm for small enough $\ep>0$.

For a collection of indices $I$ introduce a set %
$$ B_I:=\{\vx\in\vX:~~
                \forall i_1\in I~\exists i_2,\dots,i_k\in I: ~~
                   |x_{i_j}-x_{i_{j+1}}|\le\ep, ~~
                \min_{i\in I, j\notin I}|x_i-x_j|>\ep \}. $$
In words, the configurations from the set $B_I$ satisfy the
condition that all $I$-particles (i.e. those with indices from $I$)
are connected by $\ep$-chains, while all others are far enough.

The sets $\{B_I\}$ define a finite partition of $\vX$ and the map
$\Q$ differs from the identical map only on the sets $B_I$ with
$|I|>1$. Thus the signed measure $\v\nu:=\Q\vmu-\vmu$ is supported
only on the sets $B_I$ with $|I|>1$. Denote
$\t{B}_I:=\{\vx\in\vX:~~|x_i-x_j|\le\ep|I|~\forall i,j\in I \}$.
Obviously $\t{B}_I$ is a proper rectangle and
$B_I\subseteq\t{B}_I~\forall I$. Applying
Lemma~\ref{l:var-properties}(c) we get %
\bea{ \nw{\v\nu} \a=\nw{\v\nu_{|\cup_{|I|>1}B_I}}
  =\sum_{|I|>1}\nw{\v\nu_{|B_I}}
\le\sum_{|I|>1}\nw{\v\nu_{|\t{B}_I}} \\
\a\le\frac12\sum_{|I|>1} |\t{B}_I| (\ns{\vmu} + \ns{\Q^*\vmu}) \\
\a\le \frac12\sum_{|J|>1} ~(\ep|I|)^{|J|}~
    \left(1+\frac{2^N\theta}{\ga}\right)\ns{\vmu}
\toas{\ep\to0}0 .}%

The completion of the proof follows the same perturbation argument
as in the proof of Theorem~\ref{t:disc1d}. \qed

\section{Generalizations}\label{s:gen}

{\bf 1. General interaction potential}. So far we have considered
local interactions based on a somewhat non-physical model of
attraction. Indeed, normally by the attraction one means something
more close to the gravitation law. In order to include these more
general local interactions (excluded by (\ref{e:interaction}))
consider an ``interaction potential'' $U:X\to\IR$ and set
$U_\ep(x):=U(x/\ep)$. Then one defines the following
generalization of the dynamically switching interactions: %
\beq{e:interaction-gen}{(\vQ_\ep\vx)_i:=\ga x_i +
  \frac{1-\ga}{|J_i|}{\sum_{j\in J_i}x_j\cdot U_\ep(x_j-x_i)} .} %
The interaction in (\ref{e:interaction}) corresponds to the
potential $U(x)$ defined by the indicator function $1_{[-1,1]}$ of
the interval $[-1,1]$. Assuming that $U(x)\ge0$ and making some
regularity type assumptions on the potential one recovers all
results obtained in Theorems~\ref{t:main1} and \ref{t:main2}.

\bigskip\n{\bf 2. Random local dynamics}. It worth notice that the
dynamics of the local units of the multicomponent system under
consideration needs not to be deterministic. Indeed, it might be
defined by a stochastic Markov chain satisfying the weak mixing
condition (in Theorem~\ref{t:main1}) and some additional
assumptions about the induced operator acting in the space of
signed measures (in Theorem~\ref{t:main2}).

\bigskip\n{\bf 3. Infinite particle systems}. Strictly speaking
our definition of the interaction does not allow to consider
infinite particle systems since the notion of the center of
gravity is not well defined in this case. To overcome this
difficulty one may assume that there exists a certain locally
finite\footnote{i.e. the degree of each vertex is finite, but
   not necessarily uniformly bounded.} %
graph of interactions and that the dynamical switching occurs only
between neighboring elements in this graph. This means that the
sets $J_i$ satisfy the property that only elements $J$ connected
to $i$ may belong to them. All our results can be extended to this
setup.

Another and potentially more promising approach (at least in the
``independent phase'') is to consider a mean field approximation
scheme. Let $\mu$ be a probability distribution describing the
position of a particle (and assume that it is the same for all
particles). If a given particle is located at a point $x\in X$
then the mean field approximation allows to calculate the center
of gravity of the particles in the $\ep$-neighborhood $B_\ep(x)$
of this point as $\frac1{\mu(B_\ep(x))}\int_{B_\ep(x)}y d\mu(y)$.
Therefore one rewrites the interaction as %
$$ Q_{\ep,\ga,\mu} x := \ga x
         + \frac{1-\ga}{\mu(B_\ep(x))}\int_{B_\ep(x)}y d\mu(y) .$$
Denoting (for a given $\mu$) by $Q_{\ep,\ga,\mu}^*$ the induced
action of the map $Q_{\ep,\ga,\mu}$ in the space of signed
measures we obtain the description of the mean field approximation
in this space: %
\beq{e:mean}{\map_{\ep,\ga}^*\mu := Q_{\ep,\ga,\map^*\mu}^*\map^*\mu .} %
In distinction to the transfer operators considered in the
previous Sections the operator $\map_{\ep,\ga}^*$ is nonlinear
which complicates its analysis a lot.

In the simplest case when $X:=S^1$ and $\map x:={2x}$ it is easy
to show that the Lebesgue measure is $\map_{\ep,\ga}^*$-invariant
for all $\ep,\ga$ and that any measure uniformly distributed on a
periodic trajectory is invariant for small enough $\ep>0$.
Nevertheless the analysis of stability of these measures (i.e. the
construction of the analogue of the SBR measure) is a much more
delicate task. Even in this simple example one needs to develop a
special technique to study properties of the nonlinear transfer
operator. Therefore this analysis will be discussed in a separate
publication

\end{document}